\documentclass[12pt]{article}
\usepackage{amsmath,amsthm,amssymb,latexsym,mathrsfs}

\newtheorem{thm}{Theorem}[section]
\newtheorem{lem}[thm]{Lemma}
\newtheorem{cor}[thm]{Corollary}
\newtheorem{prop}[thm]{Proposition}
\newtheorem{rem}[thm]{Remark}
\newtheorem{ex}[thm]{Example}
\newtheorem*{LT}{Liggett Theorem}

\newcommand{\lrf}[1]{\lfloor #1\rfloor}
\newcommand{\la}{\lambda}
\newcommand{\A}{\mathscr{A}}
\renewcommand{\L}{\mathscr{L}}

\title
   {\bf Log-concavity and LC-positivity}
\author
   {Yi\ Wang$^{\rm a}$
   \footnote{
   {\it Email address:}\quad wangyi@dlut.edu.cn (Y. Wang)}
   {\ }\footnote{Partially supported by NSF of China 10471016}
   \quad and \quad
   Yeong-Nan\ Yeh$^{\rm b}$\footnote{Partially supported by NSC 94-2115-001-017}}
\date{\footnotesize{
   $^{\rm a}$
   Department of Applied Mathematics,
   Dalian University of Technology,
   Dalian 116024, China
   \\
   $^{\rm b}$
   Institute of Mathematics,
   Academia Sinica,
   Taipei 11529, Taiwan}}

\begin{document}

\maketitle

\begin{abstract}
A triangle $\{a(n,k)\}_{0\le k\le n}$ of nonnegative numbers is LC-positive if for each $r$,
the sequence of polynomials $\sum_{k=r}^{n}a(n,k)q^k$ is $q$-log-concave.
It is double LC-positive if both triangles $\{a(n,k)\}$ and $\{a(n,n-k)\}$ are LC-positive.
We show that if $\{a(n,k)\}$ is LC-positive
then the log-concavity of the sequence $\{x_k\}$ implies that of the sequence $\{z_n\}$ defined by
$z_n=\sum_{k=0}^{n}a(n,k)x_k$,
and if $\{a(n,k)\}$ is double LC-positive
then the log-concavity of sequences $\{x_k\}$ and $\{y_k\}$
implies that of the sequence $\{z_n\}$ defined by
$z_n=\sum_{k=0}^{n}a(n,k)x_ky_{n-k}$.
Examples of double LC-positive triangles include
the constant triangle and the Pascal triangle.
We also give a generalization of a result of Liggett
that is used to prove a conjecture of Pemantle on characteristics of negative dependence.
\end{abstract}
\section{Introduction}
\hspace*{\parindent}
Let $x_0,x_1,x_2,\ldots$ be a sequence of nonnegative numbers and with no internal zeros.
By the latter we mean that
there are no three indices $i<j<k$ such that $x_i,x_k\not=0$ and $x_j=0$.
We say that the sequence is {\it log-concave} ({\bf LC})
if $x_{i-1}x_{i+1}\le x_i^2$ for all $i>0$.
It is well known that the sequence $\{x_k\}$ is log-concave if and only if
$x_{i-1}x_{j+1}\le x_ix_j$ for all $j\ge i\ge 1$
(see \cite[Proposition 2.5.1]{Bre89} for instance),
or equivalently,
all minors of order $2$ of the infinite matrix $M=(x_{i-j})_{i,j\ge 0}$ are nonnegative
(where $x_k=0$ if $k<0$).
For this reason a log-concave sequence with no internal zeros is also called PF$_2$
(the notation actually has a precisely motivation, see~\cite{Bre89,Kar68}).
Log-concave sequences arise often in combinatorics,
algebra, geometry, analysis, probability and statistics.
There have been many attempts to develop techniques for the log-concavity problems.
We refer the reader to Stanley's survey article~\cite{Sta89}
and Brenti's supplement~\cite{Bre94} for details.

Let $\{a(n,k)\}_{0\le k\le n}$ be a triangular array of nonnegative numbers.
Define two linear transformations of sequences by
\begin{eqnarray}\label{seq1}
z_n=\sum\limits_{k=0}^na(n,k)x_k,\quad n=0,1,2,\ldots
\end{eqnarray}
and
\begin{eqnarray}\label{seq2}
z_n=\sum\limits_{k=0}^na(n,k)x_ky_{n-k},\quad n=0,1,2,\ldots
\end{eqnarray}
respectively.
We say that the linear transformation (\ref{seq1}) has {\it the PLC property}
if it preserves the log-concavity of sequences, i.e.,
the log-concavity of $\{x_n\}$ implies that of $\{z_n\}$.
We say that the linear transformation (\ref{seq2}) has {\it the double PLC property}
if the log-concavity of $\{x_n\}$ and $\{y_n\}$ implies that of $\{z_n\}$.
The corresponding triangle $\{a(n,k)\}$
is also called {\it PLC} and {\it double PLC} respectively.
Clearly, the double PLC property implies the PLC property.

It is well known that the ordinary convolution
\begin{eqnarray*}
z_n=\sum_{k=0}^nx_ky_{n-k},\quad n=0,1,2,\ldots
\end{eqnarray*}
is double PLC,
which can be obtained as a consequence of the fact
that the product of $TP_2$ matrices is $TP_2$
(see Karlin~\cite[p. 394]{Kar68} for instance)
or by a direct argument (see Menon~\cite{Men69} for instance).
Using the same fact,
Walkup can manage to prove that the binomial convolution
\begin{eqnarray*}
z_n=\sum_{k=0}^n\binom{n}{k}x_ky_{n-k},\quad n=0,1,2,\ldots
\end{eqnarray*}
is double PLC (\cite[Theorem 1]{Wal76}).
A more general result is due to Liggett
(see~\cite[Theorem 3]{Lig97} or Section 3 of this paper).
However,
there is no systematic study of linear transformations that are double PLC.
The possible reason for this is that very few examples of such linear transformations are known.
In the present paper we develop techniques to deal with the problems of
finding these kind of linear transformations
and apply these techniques to generate new log-concave sequences from existing ones.

When the triangle $\{a(n,k)\}$ is PLC,
the linear transformation (\ref{seq1})
has to send any log-concave sequence $\{x_k\}$ to a log-concave sequence $\{z_n\}$.
So, by taking the special log-concave sequence $\{x_k\}$,
we may obtain certain necessary conditions such that $\{a(n,k)\}$ is PLC
from the log-concavity of the associated sequence $\{z_n\}$.
\begin{rem}\label{rem1.1}
Let the triangle $\{a(n,k)\}$ be PLC.
Then for $r\in\mathbb{N}$ and $p>0$,
\begin{enumerate}
\item [(i)]
the column sequence $\{a(n,r)\}_{n\ge r}$ is log-concave;
\item [(ii)]
the row-sum sequence $a(n)=\sum_{k=0}^na(n,k)$ is log-concave; and
\item [(iii)]
the sequence $\A_r(n;p)=\sum_{k=r}^na(n,k)p^k$ is log-concave for $n\ge r$.
\end{enumerate}

We can view $\A_r(n;p)$ as a polynomial in $p$.
By (iii), the polynomial
$$\A^2_r(n;p)-\A_r(n-1;p)\A_r(n+1;p)$$
takes nonnegative values when $p>0$,
and so that its leading coefficient
$$a^2(n,n)-a(n-1,n-1)a(n+1,n+1)$$
has to be nonnegative.
In other words,
the diagonal sequence $\{a(n,n)\}_{n\ge 0}$ is log-concave.
\end{rem}

In order to state our sufficient conditions for $\{a(n,k)\}$ to be PLC,
we introduce some terminology and notation.
Let $q$ be an indeterminate
and $\{f_n(q)\}_{n\ge 0}$ a sequence of polynomials in $q$.
We say that the sequence $\{f_n(q)\}_{n\ge 0}$ is {\it $q$-log-concave}
if for each $n\ge 1$,
$f_n^2(q)-f_{n-1}(q)f_{n+1}(q)$ has nonnegative coefficients as a polynomial in $q$.
The concept of $q$-log-concavity is first suggested by Stanley (see \cite[p. 795]{Sag92Trans}).
We refer the reader to \cite{But90,Kra89,Ler90,Sag92DM,Sag92Trans}
for further information about $q$-log-concavity.
Now for $0\le r\le n$, define the polynomial
\begin{eqnarray*}
{\A}_r(n;q)=\sum_{k=r}^n a(n,k)q^k.
\end{eqnarray*}
We say that the triangle $\{a(n,k)\}$ has {\it the LC-positive property}
if for each $r\ge 0$,
the sequence of polynomials $\{{\A}_r(n;q)\}_{n\ge r}$ is $q$-log-concave in $n$.
(We remind the reader that the definition is different from Remark~\ref{rem1.1}~(iii).)
Define the reciprocal triangle $\{a^*(n,k)\}$ of $\{a(n.k)\}$ by
$$a^*(n,k)=a(n,n-k),\quad 0\le k\le n.$$
We say that the triangle $\{a(n,k)\}$ has {\it the double LC-positive property}
if both $\{a(n,k)\}$ and $\{a^*(n,k)\}$ have the LC-positive property.
\begin{ex}\label{ex1.1}
Consider $a(n,k)\equiv 1$ for $0\le k\le n$.
Then ${\A}_r(n;q)=\sum_{k=r}^nq^k$ for $0\le r\le n$.
It immediately follows that
$${\A}^2_r(n;q)-{\A}_r(n-1;q){\A}_r(n+1;q)=q^{n+r},$$
and so that $\{{\A}_r(n;q)\}$ is $q$-log-concave in $n$.
Thus the constant triangle $\{a(n,k)\}$ is LC-positive
and therefore double LC-positive since $a^*(n,k)=a(n,k)$.
\end{ex}
\begin{ex}\label{ex1.2}
Consider $a(n,k)=\binom{n}{k}$.
Then ${\A}_r(n;q)=\sum_{k=r}^n\binom{n}{k}q^k$.
We have
$$\A_r(n;q)
=\sum_{k=r}^{n}\left[\binom{n-1}{k}+\binom{n-1}{k-1}\right]q^k
=(q+1)\A_r(n-1;q)+\binom{n-1}{r-1}q^r.$$
It follows that
\begin{eqnarray*}
& & \A^2_r(n;q)-\A_r(n-1;q)\A_r(n+1;q)\\
&=& \A_r(n;q)\left[(q+1)\A_r(n-1;q)+\binom{n-1}{r-1}q^r\right]\\
& & -\A_r(n-1;q)\left[(q+1)\A_r(n;q)+\binom{n}{r-1}q^r\right]\\
&=&\binom{n-1}{r-1}q^r\A_r(n;q)-\binom{n}{r-1}q^r\A_r(n-1;q)\\
&=&\sum_{k=r}^{n}\left[\binom{n-1}{r-1}\binom{n}{k}-\binom{n}{r-1}\binom{n-1}{k}\right]q^{k+r}\\
&=&\sum_{k=r}^{n}\left[\binom{n-1}{r-1}\binom{n-1}{k-1}-\binom{n-1}{r-2}\binom{n-1}{k}\right]q^{k+r},
\end{eqnarray*}
which has nonnegative coefficients
by the log-concavity of the binomial coefficients.
Hence $\{\A_r(n;q)\}$ is $q$-log-concave in $n$.
Thus the Pascal triangle $\{a(n,k)\}$ is LC-positive
and therefore double LC-positive since $a^*(n,k)=a(n,k)$.
\end{ex}

The object of this paper is twofold.
First, we show that LC-positive triangles are PLC
and that double LC-positive triangles are double PLC.
Second, we present some examples of
PLC and double PLC triangles
by showing their LC-positivity.
We also give a generalization of a result of Liggett
that is used to prove a conjecture of Pemantle on characteristics of negative dependence.
\section{Theorems}
\hspace*{\parindent}
In this section we discuss the LC-positivity in detail
and establish the relation between the (double) LC-positivity and (double) PLC property.
The following simple result will be used repeatedly in our discussion.
\begin{lem}\label{abel}
Let $s\in\mathbb{P}$.
Suppose that two sequences $a_0,\ldots, a_s$ and $X_0,\ldots, X_s$ of real numbers
satisfy the following two conditions:
\begin{itemize}
\item [(a)] $\sum_{k=r}^s a_k\ge 0$ for all $0\le r\le s$;
\item [(b)] $0\le X_0\le X_1\le\ldots\le X_s$.
\end{itemize}
Then $\sum_{k=0}^sa_kX_k\ge X_0\sum_{k=0}^sa_k\ge 0$.
\end{lem}
\proof
Applying the Abel's partial summation formula
$$\sum_{k=0}^sa_kX_k
=(a_0+a_1+\cdots+a_s)X_0+(a_1+\cdots+a_s)(X_1-X_0)+\cdots+a_s(X_s-X_{s-1}),$$
the statement immediately follows.
\qed

We first consider the relation between the LC-positivity and the PLC property.
Let $\{a(n,k)\}_{0\le k\le n}$ be a triangle of nonnegative numbers
and $\{x_k\}_{k\ge 0}$ be a log-concave sequence.
It is convenient to extend the definition of $x_k$ and $a(n,k)$
by setting $x_k=0$ for $k<0$ and $a(n,k)=0$ for $k<0$ or $k>n$.
Let $\{z_n\}_{n\ge 0}$ be the sequence defined by (\ref{seq1})
and denote $\Delta_n=z_n^2-z_{n-1}z_{n+1}$.
Then we need that $\Delta_n\ge 0$ for each $n\ge 1$.
Note that
\begin{eqnarray}\label{delta}
\Delta_n
=\left\{\sum_{k=0}^n a(n,k)x_k\right\}^2
-\left\{\sum_{k=0}^{n-1} a(n-1,k)x_k\right\}
 \left\{\sum_{k=0}^{n+1} a(n+1,k)x_k\right\}
\end{eqnarray}
is a quadratic form in $n+2$ variables $x_0,x_1,\ldots,x_{n+1}$.
Such quadratic forms are generally not positive semidefinite.
Hence the log-concavity of $\{x_k\}$ is indispensable for our purposes.
To see this let us take $a(n,k)\equiv 1$ for $0\le k\le n$ as an example.
In this case we have
$$\Delta_2=(x_0+x_1)^2-x_0(x_0+x_1+x_2)=x_1^2+x_0x_1-x_0x_2.$$
Clearly,
$\Delta_2$ may take negative values for nonnegative $x_k$'s,
but must be nonnegative when $x_0,x_1,x_2$ is log-concave.

To utilize the assumption for $\{x_k\}$,
recall that $\{x_k\}$ is log-concave if and only if
$x_{i-1}x_{j+1}\le x_ix_j$ for $j\ge i\ge 1$.
In other words,
the $x_ix_j$'s with the same ``weight" $i+j$ are comparable.
Collect together those terms in $\Delta_n$ with the same weight $t$
and denote their sum by $S_t$.
For $0\le k\le \lrf{t/2}$,
let $a_k(n,t)$ be the coefficient of the term $x_kx_{t-k}$ in $\Delta_n$.
Then $\Delta_n=\sum_{t=0}^{2n}S_t$ and $S_t=\sum_{k=0}^{\lrf{t/2}}a_k(n,t)x_kx_{t-k}$.
Thus it suffices that $S_t\ge 0$ for each $0\le t\le 2n$.
Note that
$x_0x_t\le x_1x_{t-1}\le x_2x_{t-2}\le\cdots$.
Hence by Lemma \ref{abel},
it suffices that $\sum_{k=r}^{\lrf{t/2}} a_k(n,t)\ge 0$ for each $0\le r\le \lrf{t/2}$.
By (\ref{delta}),
$$a_k(n,t)=2a(n,k)a(n,t-k)-a(n-1,k)a(n+1,t-k)-a(n+1,k)a(n-1,t-k)$$
for $k<t/2$, and
$$a_k(n,t)=a^2(n,k)-a(n-1,k)a(n+1,k)$$
for $t$ even and $k=t/2$.
Denote
\begin{eqnarray}\label{arnt}
A_r(n,t)=\sum_{k=r}^{\lrf{t/2}}a_k(n,t).
\end{eqnarray}
Then it is not difficult to see that
$A_r(n,t)$ is precisely the coefficient of $q^t$
in the polynomial ${\A}^2_r(n;q)-{\A}_r(n-1;q){\A}_r(n+1;q)$,
i.e.,
\begin{eqnarray}\label{arnt1}
{\A}^2_r(n;q)-{\A}_r(n-1;q){\A}_r(n+1;q)=\sum_{t=2r}^{2n}A_r(n,t)q^t.
\end{eqnarray}
So the following lemma is immediate.
\begin{lem}\label{lc-positive}
With the notation above,
the triangle $\{a(n,k)\}_{0\le k\le n}$ is LC-positive
if and only if $A_r(n,t)\ge 0$ for all $2r\le t\le 2n$.
\end{lem}

We can  now conclude the first main result of this paper from the discussion above.
\begin{thm}\label{thm1}
The LC-positive triangles are PLC.
\end{thm}

We next relate the double LC-positivity with the double PLC property.
We need the following.
\begin{prop}\label{basic}
Given a triangle $\{a(n,k)\}_{0\le k\le n}$ of nonnegative numbers
and two log-concave sequences $\{x_k\}_{k\ge 0}$ and $\{y_k\}_{k\ge 0}$,
define three triangles $\{b(n,k)\}, \{c(n,k)\}$ and $\{d(n,k)\}$ by
$$b(n,k)=a(n,k)x_k,\quad c(n,k)=a(n,k)y_{n-k},\quad d(n,k)=a(n,k)x_ky_{n-k}.$$
For $2r\le t\le 2n$,
define $B_r(n,t),C_r(n,t)$ and $D_r(n,t)$ similar to $A_r(n,t)$ in (\ref{arnt}).
\begin{enumerate}
\item [(i)]
If the triangle $\{a(n,k)\}$ is LC-positive,
then the triangle $\{b(n,k)\}$ is LC-positive
and $B_r(n,t)\ge A_r(n,t)x_rx_{t-r}$.
\item [(ii)]
If the triangle $\{a(n,k)\}$ is double LC-positive,
then the triangle $\{c(n,k)\}$ is LC-positive
and $C_r(n,t)\ge A_r(n,t)y_{n-t+r}y_{n-r}$ for $t\le n+r$.
\item [(iii)]
If the triangle $\{a(n,k)\}$ is double LC-positive,
then the triangle $\{d(n,k)\}$ is LC-positive
and $D_r(n,t)\ge A_r(n,t)x_rx_{t-r}y_{n-t+r}y_{n-r}$ for $t\le n+r$.
\end{enumerate}
\end{prop}
\proof
Clearly, (iii) follows from (i) and (ii),
so it suffices to prove (i) and (ii).

(i)\quad
Let $0\le t\le 2n$.
It is easy to see by definition that
$b_k(n,t)=a_k(n,t)x_kx_{t-k}$ for $0\le k\le \lrf{t/2}$.
Hence for $0\le r\le \lrf{t/2}$,
\begin{eqnarray*}
B_r(n,t)=\sum_{k=r}^{\lrf{t/2}}b_k(n,t)
=\sum_{k=r}^{\lrf{t/2}}a_k(n,t)x_kx_{t-k}.
\end{eqnarray*}
Now $\{a(n,k)\}$ is LC-positive and
$x_0x_t\le x_1x_{t-1}\le x_2x_{t-2}\le\cdots$ by the log-concavity of $\{x_k\}$.
From Lemma \ref{abel} it follows that
$$B_r(n,t)\ge x_rx_{t-r}\sum_{k=r}^{\lrf{t/2}}a_k(n,t)=A_r(n,t)x_rx_{t-r}\ge 0.$$
So the triangle $\{b(n,k)\}$ is LC-positive.

(ii)\quad
Let $2r\le t\le 2n$.
We need to prove $C_r(n,t)\ge 0$.
For brevity,
we do this only for the case $t$ odd
since the same technique is still valid for the case $t$ even.

Let $t=2s+1$.
For $0\le k\le s$, denote
\begin{eqnarray*}
\alpha_k&=&a(n,k)a(n,t-k),\\
\beta_k&=&a(n-1,k)a(n+1,t-k),\\
\gamma_k&=&a(n+1,k)a(n-1,t-k),
\end{eqnarray*}
and $Y_k=y_{n-t+k}y_{n-k}$.
Then
$$a_k(n,t)=2\alpha_k-\beta_{k}-\gamma_{k}$$
and
\begin{eqnarray*}
c_k(n,t)=2\alpha_kY_k-\beta_{k}Y_{k+1}-\gamma_{k}Y_{k-1}
\end{eqnarray*}
by definition.
It follows that
\begin{eqnarray*}
C_r(n,t)
& = & \sum_{k=r}^{s}(2\alpha_kY_k-\beta_{k}Y_{k+1}-\gamma_{k}Y_{k-1})\\
& = & \sum_{k=r}^{s}(2\alpha_k-\beta_{k-1}-\gamma_{k+1})Y_k+\beta_{r-1}Y_r-\gamma_rY_{r-1},
\end{eqnarray*}
where we use the fact that $Y_{s+1}=Y_s$ and $\gamma_{s+1}=\beta_s$.
Note that $\{Y_k\}$ is nondecreasing by the log-concavity of $\{y_k\}$ and
\begin{eqnarray*}
\lefteqn{ 2\alpha_k-\beta_{k-1}-\gamma_{k+1}
= 2a^*(n,n-k)a^*(n,n-t+k)}\\
& & -a^*(n-1,n-k)a^*(n+1,n-t+k)
 -a^*(n+1,n-k)a^*(n-1,n-t+k)\\
&=& a^*_{n-t+k}(n,2n-t).
\end{eqnarray*}
Hence by the LC-positivity of $\{a^*(n,k)\}$, we have
\begin{eqnarray}\label{crnt}
C_r(n,t)
& = & \sum_{j=n-t+r}^{\lrf{(2n-t)/2}}a^*_{j}(n,2n-t)Y_{j-n+t}+\beta_{r-1}Y_r-\gamma_rY_{r-1}\nonumber\\
& \ge & Y_r\sum_{j=n-t+r}^{\lrf{(2n-t)/2}}a^*_{j}(n,2n-t)+\beta_{r-1}Y_r-\gamma_rY_{r-1}\nonumber\\
& = & Y_r\sum_{k=r}^{s}(2\alpha_k-\beta_{k-1}-\gamma_{k+1})+\beta_{r-1}Y_r-\gamma_rY_{r-1}\nonumber\\
& = & Y_r\sum_{k=r}^{s}(2\alpha_k-\beta_{k}-\gamma_{k})+\gamma_r(Y_r-Y_{r-1})\nonumber\\
& = & A_r(n,t)Y_r+\gamma_r(Y_r-Y_{r-1}).
\end{eqnarray}
Thus $C_r(n,t)\ge A_r(n,t)y_{n-t+r}y_{n-r}\ge 0$ since $Y_r\ge Y_{r-1}$, as desired.
\qed

Now we present the second main result of this paper.
\begin{thm}\label{thm2}
The double LC-positive triangles are double PLC.
\end{thm}
\proof
Let the triangle $\{a(n,k)\}$ be double LC-positive.
Suppose that both $\{x_k\}$ and $\{y_k\}$ are log-concave.
Then the triangle $\{a(n,k)x_ky_{n-k}\}$ is LC-positive by Proposition~\ref{basic} (iii)
and is therefore PLC by Theorem~\ref{thm1}.
Thus the row-sum sequence
$$z_n=\sum_{k=0}^na(n,k)x_ky_{n-k},\quad n=0,1,2,\ldots$$
is log-concave.
In other words,
the triangle $\{a(n,k)\}$ is double PLC.
\qed

We can give some more practicable conditions
that imply the LC-positivity.
We have seen that Lemma~\ref{abel},
especially Condition~(a),
plays a key role in the proof of the LC-positivity of Proposition~\ref{basic}.
Clearly,
Condition~(a) is implied by the following two conditions:
\begin{enumerate}
\item [(a1)]
$a_0,a_1,\ldots,a_s$ changes from nonpositive to nonnegative values;
\item [(a2)]
$\sum_{k=0}^s a_k\ge 0$.
\end{enumerate}

These two conditions are easier to check than Condition~(a).
For example, Condition~(a1) can be obtained by showing that
the sequence $\{a_k\}$ is nondecreasing and eventually nonnegative.
In this case the analytic tools are often effective.
On the other hand,
Condition~(a2) is just the simplest one of inequalities in Condition~(a)
and the methods of generating functions will be useful
(see~\cite{Wil94} for details).
By Lemma~\ref{lc-positive},
$\{a(n,k)\}$ is LC-positive if and only if the inequality
$\sum_{k=r}^{\lrf{t/2}} a_k(n,t)\ge 0$
for all $2r\le t\le 2n$,
so the following corollary is immediate.
\begin{cor}\label{cor2.1}
Suppose that the following two conditions hold:
\begin{enumerate}
\item [(A)]
There exists an index $m=m(n,t)$ such that
$a_k(n,t)<0$ for $k<m$ and $a_k(n,t)\ge 0$ for $k\ge m$;
\item [(B)]
The sequence $\{\A_0(n;q)\}_{n\ge 0}$ is $q$-log-concave.
\end{enumerate}
Then the triangle $\{a(n,k)\}$ is LC-positive and therefore PLC.
\end{cor}
\begin{cor}
Suppose that the triangle $\{a(n,k)\}$ satisfies Condition~(A) and (B)
in Corollary~\ref{cor2.1}
and $\{a^*(n,k)\}$ satisfies Condition~(A).
Then $\{a(n,k)\}$ is double LC-positive and therefore double PLC.
\end{cor}
\proof
Clearly, it suffices to show that $\{\A^*_0(n;q)\}$ is $q$-log-concave.
We have
\begin{eqnarray*}
\A^*_0(n;q)=\sum_{k=0}^na(n,n-k)q^k
=\sum_{k=0}^na(n,k)q^{n-k}=q^n\A_0(n;q^{-1}).
\end{eqnarray*}
It follows that
\begin{eqnarray*}
& & {\A^*_0}^2(n;q)-\A^*_0(n-1;q)\A^*_0(n+1;q)\\
&=& q^{2n}\left[\A^2_0(n;q^{-1})-\A_0(n-1;q^{-1})\A_0(n+1;q^{-1})\right],
\end{eqnarray*}
which has nonnegative coefficients by the $q$-log-concavity of $\{\A_0(n;q)\}$,
as desired.
\qed
\section{Applications}
\hspace*{\parindent}
In this section
we give some examples of PLC and double PLC triangles
by showing their LC-positivity.
In particular,
we give a generalization of a result of Liggett
that is used to prove a conjecture of Pemantle on characteristics of negative dependence.

Denote by $\mathfrak{S}$ the set of sequences $\{u_k\}_{k\in\mathbb{Z}}$ of nonnegative numbers.
Given two nonnegative numbers $\la$ and $\mu$,
define the linear operator $\L=\L[\la,\mu]$ on $\mathfrak{S}$ by
$$\L(u_k)=\la u_k+\mu u_{k-1},\quad k\in\mathbb{Z}.$$
For $n\ge 2$, define $\L^n=\L(\L^{n-1})$ by induction.
It is convenient to view $\L^0$ as the identity operator.
Let the sequence $\{u_k\}_{k\in\mathbb{Z}}$ be log-concave.
Then the sequence $\{\L(u_k)\}_{k\in\mathbb{Z}}$ is also log-concave since
\begin{eqnarray*}
& & [\L(u_k)]^2-\L(u_{k-1})\L(u_{k+1})\\
&=& (\la u_k+\mu u_{k-1})^2-(\la u_{k-1}+\mu u_{k-2})(\la u_{k+1}+\mu u_{k})\\
&=& \la^2(u_k^2-u_{k-1}u_{k+1})+\la\mu (u_{k-1}u_k-u_{k-2}u_{k+1})+\mu^2(u_{k-1}^2-u_{k-2}u_k).
\end{eqnarray*}
Thus we can conclude by induction
that the sequence $\{\L^n(u_k)\}_{k\in\mathbb{Z}}$ is log-concave for each $n\ge 0$.
\begin{thm}\label{thm-back}
Given two nonnegative numbers $\la,\mu$
and a log-concave sequence $\{u_k\}$,
define $a(n,k)=\L^n[\la,\mu](u_k)$ for $0\le k\le n$.
Then the triangle $\{a(n,k)\}_{0\le k\le n}$ is double LC-positive and therefore double PLC.
\end{thm}
\proof
Denote $a_k=\L^{n-1}[\la,\mu](u_k)$ for $k\in\mathbb{Z}$.
Then the sequence $\{a_k\}_{k\in\mathbb{Z}}$ is log-concave and
${\A}_r(n-1;q)=\sum_{k=r}^{n-1} a_kq^k$.
We have
\begin{eqnarray*}
{\A}_r(n;q)
&=& \sum_{k=r}^{n} (\la a_k+\mu a_{k-1})q^k\\
&=& \la\sum_{k=r}^{n}a_kq^k+\mu\sum_{k=r}^{n}a_{k-1}q^k\\
&=& (\la +\mu q){\A}_r(n-1;q)+\la a_{n}q^n+\mu a_{r-1}q^r,
\end{eqnarray*}
and similarly,
\begin{eqnarray*}
{\A}_r(n+1;q)=(\la +\mu q){\A}_r(n;q)+\la(\la a_{n+1}+\mu a_n)q^{n+1}+\mu(\la a_{r-1}+\mu a_{r-2})q^r.
\end{eqnarray*}
It follows that
\begin{eqnarray}\label{arnt-formula}
& &  {\A}^2_r(n;q)-{\A}_r(n-1;q){\A}_r(n+1;q)\nonumber\\
&=&  \A_r(n;q)\left[(\la +\mu q){\A}_r(n-1;q)+\la a_{n}q^n+\mu a_{r-1}q^r\right]\nonumber\\
& & -{\A}_r(n-1;q)\left[(\la +\mu q){\A}_r(n;q)
    +\la(\la a_{n+1}+\mu a_n)q^{n+1}+\mu(\la a_{r-1}+\mu a_{r-2})q^r\right]\nonumber\\
&=& (\la a_{n}q^n+\mu a_{r-1}q^r){\A}_r(n;q)\nonumber\\
& & -\left[\la(\la a_{n+1}+\mu a_n)q^{n+1}
    +\mu(\la a_{r-1}+\mu a_{r-2})q^r\right]{\A}_r(n-1;q)\nonumber\\
&=&  \la\sum_{k=r}^{n}(\la a_k+\mu a_{k-1})a_{n}q^{n+k}
    +\mu\sum_{k=r}^{n}a_{r-1}(\la a_k+\mu a_{k-1})q^{k+r}\nonumber\\
& & -\la\sum_{k=r}^{n-1}a_k(\la a_{n+1}+\mu a_n)q^{n+k+1}
    -\mu\sum_{k=r}^{n-1}(\la a_{r-1}+\mu a_{r-2})a_kq^{k+r}\nonumber\\
&=&  \la^2\sum_{k=r+1}^{n}(a_{k}a_{n}-a_{k-1}a_{n+1})q^{n+k}
    +\mu^2\sum_{k=r}^{n}(a_{r-1}a_{k-1}-a_{r-2}a_{k})q^{k+r}\nonumber\\
& & +(\la^2 a_r+2\la\mu a_{r-1}+\mu^2 a_{r-2})a_{n}q^{n+r},
\end{eqnarray}
which  has nonnegative coefficients by the log-concavity of $\{a_k\}$.
Hence the triangle $\{a(n,k)\}_{0\le k\le n}$ is LC-positive.

On the other hand,
let $u^*_k=u_{-k}$ for $k\in\mathbb{Z}$.
Then the sequence $\{u^*_k\}_{k\in\mathbb{Z}}$ is log-concave
and $a^*(n,k)=\L^n[\mu,\la](u^*_k)$.
Thus the triangle $\{a^*(n,k)\}_{0\le k\le n}$ is also LC-positive,
and the triangle $\{a(n,k)\}_{0\le k\le n}$ is therefore double LC-positive.
\qed
\begin{rem}\label{rem-arnq}
Let the triangle $\{a(n,k)\}$ be the same as Theorem \ref{thm-back}.
Then by (\ref{arnt1}) and (\ref{arnt-formula}), the inequality
$$A_r(n,t)\ge (a_{r-1}a_{t-r-1}-a_{r-2}a_{t-r})\mu^2$$
holds for $t\le n+r$
(the equality holds when $t<n+r$).
We will use this inequality repeatedly in the proof of Theorem~\ref{lig-gen}.
\end{rem}

Taking $\la=\mu=1/2$ and $u_k\equiv 1$ in Theorem~\ref{thm-back}
leads to the following well-known result.
\begin{cor}\label{ordconv}
If the sequences $\{x_n\}$ and $\{y_n\}$ are log-concave,
then so is their ordinary convolution
$z_n=\sum_{k=0}^nx_ky_{n-k},\quad n=0,1,2,\ldots$.
\end{cor}
\begin{cor}\label{wylaa}
Let $a,b$ be two nonnegative integers and $a\ge b$.
If the sequences $\{x_n\}$ and $\{y_n\}$ are log-concave,
then so is the sequence
$$z_n=\sum_{k=0}^n\binom{a+n}{b+k}x_ky_{n-k},\quad n=0,1,2,\ldots.$$
\end{cor}
\proof
The statement follows by taking $\la=\mu=1$ and $u_k=\binom{a}{b+k}$ in Theorem~\ref{thm-back}.
(We remind the reader that $\binom{n}{k}=0$ unless $0\le k\le n$.)
\qed

A special interesting case of Corollary~\ref{wylaa} is the following.
\begin{cor}\label{expconv}
If the sequences $\{x_n\}$ and $\{y_n\}$ are log-concave,
then so is their binomial convolution
$z_n=\sum_{k=0}^n\binom{n}{k}x_ky_{n-k},\quad n=0,1,2,\ldots$.
\end{cor}
\begin{rem}
Corollary~\ref{ordconv}, \ref{wylaa} and \ref{expconv}
can also be followed directly from Theorem~\ref{thm2}
by showing the double LC-positivity of the associated triangles.
Actually, the double LC-positivity of the constant triangle and the Pascal triangle
have been shown in Example~\ref{ex1.1} and \ref{ex1.2} respectively.
In \cite{WYlaa},
we showed the LC-positivity of the triangle $a(n,k)=\binom{a+n}{b+k}$ for $0\le k\le n$
by showing that Condition~(A) and (B) in Corollary~\ref{cor2.1} are satisfied.
This result can also be followed
by the same technique used in Example~\ref{ex1.2}.
Note that
$$a^*(n,k)=\binom{a+n}{b+(n-k)}=\binom{a+n}{(a-b)+k}.$$
Hence $a^*(n,k)$ is also LC-positive.
Thus the triangle $\{a(n,k)\}$ is double LC-positive.
\end{rem}

It is easy to extend Corollary~\ref{expconv} by induction
to several log-concave sequences.
\begin{cor}
If $\ell$ sequences $\{x_k^{(1)}\},\{x_k^{(2)}\},\ldots,\{x_k^{(\ell)}\}$ are all log-concave,
then so is the sequence
$$X_n=\sum\binom{n}{k_1,k_2,\ldots,k_{\ell}}
x_{k_1}^{(1)}x_{k_2}^{(2)}\cdots x_{k_{\ell}}^{(\ell)},\quad n=0,1,2,\ldots,$$
where the sum is over all nonnegative integers $k_1,\ldots,k_{\ell}$ such that $k_1+k_2\cdots+k_{\ell}=n$.
\end{cor}

The following theorem is in a sense ``dual" to Theorem~\ref{thm-back}.
\begin{thm}\label{thm-for}
Let $\alpha,\beta$ be two nonnegative numbers
and $\{a(n,k)\}_{0\le k\le n}$ a triangle of nonnegative numbers.
Suppose that each row of $\{a(n,k)\}$ is log-concave and satisfies the recurrence relation
\begin{eqnarray}\label{formula-for}
a(n,k)=\alpha a(n+1,k)+\beta a(n+1,k+1),\quad k=0,1,\ldots,n.
\end{eqnarray}
Then the triangle $\{a(n,k)\}$ is double LC-positive and therefore double PLC.
\end{thm}
\proof
Denote $a(n+1,k)=v_k$ for $0\le k\le n+1$.
Then the sequence $\{v_k\}$ is log-concave
and ${\A}_r(n+1;q)=\sum_{k=r}^{n+1}v_kq^k$.
By the recurrence relation (\ref{formula-for}) we have
$${\A}_r(n;q)
=\sum_{k=r}^{n}(\alpha v_k+\beta v_{k+1})q^k
=(\alpha +\beta q^{-1}){\A}_r(n+1;q)-\alpha v_{n+1}q^{n+1}-\beta v_rq^{r-1},$$
and similarly,
\begin{eqnarray*}
{\A}_r(n-1;q)=(\alpha +\beta q^{-1}){\A}_r(n;q)
-\alpha (\alpha v_{n}+\beta v_{n+1})q^{n}-\beta (\alpha v_r+\beta v_{r+1})q^{r-1}.
\end{eqnarray*}
It follows that
\begin{eqnarray*}
& & {\A}^2_r(n;q)-{\A}_r(n+1;q){\A}_r(n-1;q)\\
&=& {\A}_r(n;q)\left[(\alpha +\beta q^{-1}){\A}_r(n+1;q)-\alpha v_{n+1}q^{n+1}-\beta v_rq^{r-1}\right]\\
& & -{\A}_r(n+1;q)\left[(\alpha +\beta q^{-1}){\A}_r(n;q)
   -\alpha (\alpha v_{n}+\beta v_{n+1})q^{n}-\beta (\alpha v_r+\beta v_{r+1})q^{r-1}\right]\\
&=& \left[\alpha (\alpha v_{n}+\beta v_{n+1})q^{n}
    +\beta (\alpha v_r+\beta v_{r+1})q^{r-1}\right]{\A}_r(n+1;q)\\
& & -(\alpha v_{n+1}q^{n+1}+\beta v_rq^{r-1}){\A}_r(n;q)\\
&=& \alpha^2\sum_{k=r+1}^{n}(v_kv_n-v_{k-1}v_{n+1})q^{n+k}
   +\beta^2\sum_{k=r}^{n}(v_{r+1}v_{k+1}-v_rv_{k+2})q^{r+k}\\
& &+v_r(\alpha^2v_n+2\alpha\beta v_{n+1}+\beta^2v_{n+2})q^{n+r},
\end{eqnarray*}
which has nonnegative coefficients by the log-concavity of $\{v_k\}$.
So the triangle $\{a(n,k)\}$ is LC-positive.

Clearly,
the reciprocal triangle $\{a^*(n,k)\}$
possesses the same property as $\{a(n,k)\}$ does.
Hence $\{a^*(n,k)\}$ is also LC-positive.
Thus the triangle $\{a(n,k)\}$ is double LC-positive.
\qed

In Theorem~\ref{thm-for},
taking $\alpha=\beta=1/2$ and $a(n,k)\equiv 1$ for $0\le k\le n$
leads to Corollary~\ref{ordconv};
and taking $\alpha=\beta=1$ and $a(n,k)=\binom{a-n}{b-k}$ for $0\le k\le n$
leads to the following.
\begin{cor}\label{--}
Let $a,b\in\mathbb{N}$ and $a\ge b$.
If the sequences $\{x_k\}$ and $\{y_k\}$ are log-concave,
then so is the sequence
$$z_n=\sum_{k=0}^n\binom{a-n}{b-k}x_ky_{n-k},\quad n=0,1,2,\ldots.$$
\end{cor}

In what follows we generalize a result of Liggett.
Let $\{x_k\}_{k\ge 0}$ be a sequence of nonnegative numbers and with no internal zeros.
Following Pemantle~\cite{Pem00} and Liggett~\cite{Lig97},
the sequence is {\it ultra-log-concave of order $m$} ({\bf ULC($m$)})
if $x_k=0$ for $k>m$ and the sequence
$\left\{x_k/\binom{m}{k}\right\}_{k=0}^m$
is log-concave.
The sequence $\{x_k\}_{k\ge 0}$ is {\it ULC($\infty$)}
if the sequence $\{k!x_k\}_{k\ge 0}$ is log-concave.
It is clear from definitions that ULC($m$) implies ULC($\ell$) for $0\le m\le \ell\le \infty$.
The concept of ultra-log-concavity is closely related to
negatively dependent Bernoulli sequences
(see \cite{Pem00} for details).
Pemantle speculates that ultra-log-concavity is characteristic of
negative dependence in the exchangeable case.
This leads to a conjecture that
the ordinary convolution of a ULC($m$) sequence and a ULC($\ell$) sequence is ULC($m+\ell$)
where $m$ and $\ell$ may be infinity
(\cite[Conjecture 7]{Pem00}).
It is not difficult to see that
the conjecture actually consists of two parts:
\begin{enumerate}
\item [(i)]
The Pascal triangle $\left\{\binom{n}{k}\right\}$ is double PLC;
\item [(ii)]
The triangle $\left\{\binom{n}{k}\binom{a-n}{b-k}\right\}$ is double PLC.
\end{enumerate}

Liggett verified the conjecture
by establishing the following stronger result.
\begin{LT}[\cite{Lig97}]
Given three log-concave sequences $\{v_k\}$, $\{x_k\}$ and $\{y_k\}$,
let
\begin{eqnarray*}
z_{n-1}&=&\sum_{k=0}^{n-1}\binom{n-1}{k}(v_k+2v_{k+1}+v_{k+2})x_ky_{n-1-k},\\
z_n&=&\sum_{k=0}^{n}\binom{n}{k}(v_k+v_{k+1})x_ky_{n-k},\\
z_{n+1}&=&\sum_{k=0}^{n+1}\binom{n+1}{k}v_kx_ky_{n+1-k}.
\end{eqnarray*}
Then $z_{n-1}z_{n+1}\le z_n^2$.
\end{LT}

Liggett's proof for his theorem,
essentially using the double LC-positivity of the Pascal triangle,
is not simple.
To see his idea more clearly, we show the following more general result.
\begin{thm}\label{lig-gen}
Given four nonnegative numbers $\alpha,\beta,\la,\mu$
and four log-concave sequences
$\{u_k\}_{k\in\mathbb{Z}}$, $\{v_k\}_{k\ge 0}$, $\{x_k\}_{k\ge 0}$ and $\{y_k\}_{k\ge 0}$,
let $a(n,k)=\L^n[\la,\mu](u_k)$ and
\begin{eqnarray*}
z_{n-1}&=&\sum_{k=0}^{n-1}a(n-1,k)(\alpha^2 v_k+2\alpha\beta v_{k+1}+\beta^2 v_{k+2})x_ky_{n-1-k},\\
z_n&=&\sum_{k=0}^{n}a(n,k)(\alpha v_k+\beta v_{k+1})x_ky_{n-k},\\
z_{n+1}&=&\sum_{k=0}^{n+1}a(n+1,k)v_kx_ky_{n+1-k}.
\end{eqnarray*}
Then $z_{n-1}z_{n+1}\le z_n^2$.
\end{thm}
\proof
Clearly,
$z_n^2-z_{n-1}z_{n+1}$ can be viewed as a quadratic form
in $n+2$ variables $v_0,v_1,\ldots,v_{n+1}$.
Let
$$z_n^2-z_{n-1}z_{n+1}=\sum_{t=0}^{2(n+1)}\sum_{k=0}^{\lrf{t/2}}e_k(n,t)v_kv_{t-k}.$$
Then we need to show that
$\sum_{k=r}^{\lrf{t/2}}e_k(n,t)\ge 0$
for $2r\le t\le 2(n+1)$.
For brevity, we do this only for the case $t$ odd.
Let $t=2s+1$.

Define $d(n,k)=a(n,k)x_ky_{n-k}$ for $0\le k\le n$.
For convenience,
set $x_k=y_k=0$ for $k<0$ and $d(n,k)=0$ for $k<0$ or $k>n$.
The triangle $\{a(n,k)\}$ is double LC-positive by Theorem~\ref{thm-back},
and so is the triangle $\{d(n,k)\}$ by Proposition~\ref{basic}.
Rewrite
\begin{eqnarray*}
z_{n-1}&=&\sum_{k=0}^{n+1}[\alpha^2 d(n-1,k)+2\alpha\beta d(n-1,k-1)+\beta^2 d(n-1,k-2)]v_k,\\
z_{n}&=&\sum_{k=0}^{n+1}[\alpha d(n,k)+\beta d(n,k-1)]v_k,\\
z_{n+1}&=&\sum_{k=0}^{n+1}d(n+1,k)v_k.
\end{eqnarray*}
Then
\begin{eqnarray*}
\lefteqn{ e_k(n,t)
= 2[\alpha d(n,k)+\beta d(n,k-1)][\alpha d(n,t-k)+\beta d(n,t-k-1)] }\\
& & -[\alpha^2 d(n-1,k)+2\alpha\beta d(n-1,k-1)+\beta^2 d(n-1,k-2)]d(n+1,t-k)\\
& & -d(n+1,k)[\alpha^2 d(n-1,t-k)+2\alpha\beta d(n-1,t-k-1)+\beta^2 d(n-1,t-k-2)]\\
&=& \alpha^2 P_k+2\alpha\beta Q_k+\beta^2 R_k,
\end{eqnarray*}
where
\begin{eqnarray*}
P_k&=&2d(n,k)d(n,t-k)-d(n-1,k)d(n+1,t-k)-d(n+1,k)d(n-1,t-k),\\
Q_k&=&d(n,k)d(n,t-k-1)+d(n,k-1)d(n,t-k)\\
   & &-d(n-1,k-1)d(n+1,t-k)-d(n+1,k)d(n-1,t-k-1),\\
R_k&=&2d(n,k-1)d(n,t-k-1)\\
   & &-d(n-1,k-2)d(n+1,t-k)-d(n+1,k)d(n-1,t-k-2).
\end{eqnarray*}
Thus it suffices to show the inequality
\begin{eqnarray}\label{e>0}
\alpha^2 \sum_{k=r}^{s}P_k
+2\alpha\beta \sum_{k=r}^{s}Q_k
+\beta^2 \sum_{k=r}^{s}R_k
\ge 0.
\end{eqnarray}

Note that $P_k=d_k(n,t)$ and $R_k=d^*_{n-t+k+1}(n,2n-t+2)$.
Hence both
\begin{eqnarray}\label{p-formula}
\sum_{k=r}^{s}P_k=D_r(n,t)
\end{eqnarray}
and
\begin{eqnarray}\label{r-formula}
\sum_{k=r}^{s}R_k=D^*_{n-t+r+1}(n,2n-t+2)
\end{eqnarray}
are nonnegative
by the double LC-positivity of the triangle $\{d(n,k)\}$.
Also,
\begin{eqnarray}\label{q-formula}
\lefteqn{ \sum_{k=r}^{s}Q_k
= \sum_{k=r}^{s}[d(n,k)d(n,t-k-1)+d(n,k-1)d(n,t-k) }\nonumber\\
& & -d(n-1,k-1)d(n+1,t-k)-d(n+1,k)d(n-1,t-k-1)]\nonumber\\
&=& [d^2(n,s)-d(n-1,s)d(n+1,s)]+\sum_{k=r-1}^{s-1}[2d(n,k)d(n,t-1-k)\nonumber\\
& & -d(n-1,k)d(n+1,t-1-k)-d(n+1,k)d(n-1,t-1-k)]\nonumber\\
& & +[d(n+1,r-1)d(n-1,t-r)-d(n,r-1)d(n,t-r)]\nonumber\\
&=& D_{r-1}(n,t-1)+[d(n+1,r-1)d(n-1,t-r)-d(n,r-1)d(n,t-r)].
\end{eqnarray}

Assume that $r=0$ or $t>n+r$.
Then $\sum_{k=r}^{s}Q_k=D_{r-1}(n,t-1)\ge 0$.
Thus the inequality~(\ref{e>0}) is trivial.
So let $r\ge 1$ and $t\le n+r$.

If we can show that there exists a nonnegative number $E=E(n,t,r)$ such that
\begin{eqnarray}\label{pqr}
\begin{cases}
\sum_{k=r}^{s}P_k \ge \mu^2 Ex_{r}x_{t-r}y_{n-t+r}y_{n-r},\\
\sum_{k=r}^{s}Q_k \ge -\la\mu Ex_{r-1}x_{t-r}y_{n-t+r}y_{n-r+1},\\
\sum_{k=r}^{s}R_k \ge \la^2 Ex_{r-1}x_{t-r-1}y_{n-t+r+1}y_{n-r+1},
\end{cases}
\end{eqnarray}
then the arithmetic-geometric mean inequality
and the log-concavity of $\{x_k\}$ and $\{y_k\}$ will give
$$\alpha^2 \sum_{k=r}^{s}P_k+\beta^2 \sum_{k=r}^{s}R_k\ge -2\alpha\beta \sum_{k=r}^{s}Q_k,$$
the required inequality.
So, to prove (\ref{e>0}), it suffices to prove (\ref{pqr}).

We use Proposition~\ref{basic} to estimate the lower bounds for
$\sum_{k=r}^{s}P_k$, $\sum_{k=r}^{s}Q_k$ and $\sum_{k=r}^{s}R_k$.
From (\ref{p-formula}) and Proposition~\ref{basic}~(iii) it is immediate that
\begin{eqnarray}\label{pk}
\sum_{k=r}^{s}P_k\ge A_r(n,t)x_{r}x_{t-r}y_{n-t+r}y_{n-r}.
\end{eqnarray}

Similarly, note that $d^*(n,k)=a^*(n,k)y_kx_{n-k}$,
it follows from (\ref{r-formula}) and Proposition~\ref{basic}~(iii) that
\begin{eqnarray}\label{rk}
\sum_{k=r}^{s}R_k
\ge  A^*_{n-t+r+1}(n,2n-t+2)y_{n-t+r+1}y_{n-r+1}x_{t-r-1}x_{r-1}.
\end{eqnarray}

To get an analogous lower bound for $\sum_{k=r}^sQ_k$ from (\ref{q-formula}),
let $c(n,k)=a(n,k)y_{n-k}$.
Then $d(n,k)=c(n,k)x_k$ and so
$$D_{r-1}(n,t-1)\ge  C_{r-1}(n,t-1)x_{r-1}x_{t-r}$$
by Proposition~\ref{basic}~(i).
However,
\begin{eqnarray*}
C_{r-1}(n,t-1)
& \ge & A_{r-1}(n,t-1)y_{n-t+r}y_{n-r+1}\\
& & +a(n+1,r-1)a(n-1,t-r)(y_{n-t+r}y_{n-r+1}-y_{n-t+r-1}y_{n-r+2})
\end{eqnarray*}
by the inequality (\ref{crnt}).
Hence we have by (\ref{q-formula})
\begin{eqnarray}\label{qk}
\sum_{k=r}^{s}Q_k
&\ge & [A_{r-1}(n,t-1)y_{n-t+r}y_{n-r+1}\nonumber\\
& & +a(n+1,r-1)a(n-1,t-r)(y_{n-t+r}y_{n-r+1}-y_{n-t+r-1}y_{n-r+2})]x_{r-1}x_{t-r}\nonumber\\
& & +[a(n+1,r-1)x_{r-1}y_{n-r+2}a(n-1,t-r)x_{t-r}y_{n-t+r-1}\nonumber\\
& & -a(n,r-1)x_{r-1}y_{n-r+1}a(n,t-r)x_{t-r}y_{n-t+r}]\nonumber\\
&=& Qx_{r-1}x_{t-r}y_{n-t+r}y_{n-r+1},
\end{eqnarray}
where
\begin{eqnarray}\label{q-formular}
Q=A_{r-1}(n,t-1)+a(n+1,r-1)a(n-1,t-r)-a(n,r-1)a(n,t-r).
\end{eqnarray}

It remains to show that three coefficients
$A_r(n,t)$, $A^*_{n-t+r+1}(n,2n-t+2)$ and $Q$
in inequalities~(\ref{pk}), (\ref{rk}) and (\ref{qk})
have the lower bounds of the forms in (\ref{pqr}).
We do this by Remark~\ref{rem-arnq}.

Denote $a_k=\L^{n-1}[\la,\mu](u_k)$.
It follows from Remark~\ref{rem-arnq} that
$$A_r(n,t)\ge (a_{r-1}a_{t-r-1}-a_{r-2}a_{t-r})\mu^2$$
and that
\begin{eqnarray*}
Q
&\ge & (a_{r-2}a_{t-r-1}-a_{r-3}a_{t-r})\mu^2
 +(\la^2 a_{r-1}+2\la\mu a_{r-2}+\mu^2 a_{r-3})a_{t-r}\\
& & -(\la a_{r-1}+\mu a_{r-2})(\la a_{t-r}+\mu a_{t-r-1})\\
&=& -(a_{r-1}a_{t-r-1}-a_{r-2}a_{t-r})\la\mu
\end{eqnarray*}
by (\ref{q-formular}).
Also, note that $a^*(n,k)=\L^n[\mu,\la](u_{-k})$.
Again by Remark~\ref{rem-arnq},
\begin{eqnarray*}
A^*_{n-t+r+1}(n,2n-t+2)
&\ge & [a^*(n-1,n-t+r)a^*(n-1,n-r)\\
& &   -a^*(n-1,n-t+r-1)a^*(n-1,n-r+1)]\la^2\\
&=&(a_{r-1}a_{t-r-1}-a_{r-2}a_{t-r})\la^2.
\end{eqnarray*}

Finally, recall that the sequence $\{a_k\}_{k\in\mathbb{Z}}$ is log-concave,
so for $r\le \lrf{t/2}$,
$$E=a_{r-1}a_{t-r-1}-a_{r-2}a_{t-r}\ge 0,$$
as required.
This completes our proof.
\qed
\begin{rem}
Let $\{a(n,k)\}$ and $\{a'(n,k)\}$ be the double LC-positive triangles
appearing in Theorem~\ref{thm-back} and \ref{thm-for} respectively.
Although the triangle $\{a(n,k)a'(n,k)\}$ is not double LC-positive in general,
it is double PLC by Theorem~\ref{lig-gen}.
\end{rem}
\section{Concluding remarks}
\hspace*{\parindent}
In this paper we provide some sufficient conditions
for linear and bilinear transformations preserving the log-concavity.
As shown in Remark~\ref{rem1.1}~(iii),
the LC-positivity is ``almost" necessary for the PLC property.
It is a challenge to give a necessary and sufficient condition for the PLC property.
On the other hand,
we believe that the techniques developed in the present paper
can be used to deal with various combinatorial inequalities.
For example,
it is possible that the log-convexity problems can be treated with the same approach.
\section*{Acknowledgements}
\hspace*{\parindent}
This work was completed during Y. Wang's stay in the Institute of Mathematics, Academia Sinica, Taipei.
He would like to thank the Institute for its support.

The authors thank the anonymous referees for their helpful comments.

\end{document}